\documentclass[12pt]{amsart}

\usepackage{amsmath}
\usepackage{amssymb}  % defines \nmid, for example
\usepackage{latexsym} % for \Box
\usepackage{comment}
\usepackage{url}
\usepackage{fullpage,url,amssymb,amsmath,amsthm,amsfonts,mathrsfs}
\usepackage[usenames,dvipsnames]{color}
\usepackage[pagebackref = true, colorlinks = true, linkcolor = blue, citecolor = Green]{hyperref}
\usepackage[alphabetic,lite]{amsrefs}
\usepackage{enumitem}
\usepackage{amscd}   % for commutative diagrams
\usepackage[all, cmtip]{xy} % for complicated commutative diagrams
\usepackage{xfrac}

\newcommand{\Q}{{\mathbb Q}}
\newcommand{\surface}{{\mathcal S}}
\newcommand{\curveC}{{\mathcal C}}
\newcommand{\curveD}{{\mathcal D}}

\numberwithin{equation}{section}

\begin{document}

\title{Zero-cycles of degree one on Skorobogatov's bielliptic surface}

\author{Brendan Creutz}
\address{School of Mathematics and Statistics, University of Canterbury, Private Bag 4800, Christchurch 8140, New Zealand}
\email{brendan.creutz@canterbury.ac.nz}
\urladdr{http://www.math.canterbury.ac.nz/\~{}bcreutz}

\begin{abstract}
Skorobogatov constructed a bielliptic surface which is a counterexample to the Hasse principle not explained by the Brauer-Manin obstruction. We show that this surface has a $0$-cycle of degree $1$, as predicted by a conjecture of Colliot-Th\'el\`ene.
\end{abstract}

\maketitle

Consider the smooth projective surface $\surface/\Q$ given by the affine equations,
	\[
		\surface : (x^2 + 1)y^2 = (x^2 + 2)z^2 = 3(t^4 - 54t^2 -117t -243)\;.
	\]
	Skorobogatov \cite{Skorobogatov99} showed that the set $\surface(\Q)$ of rational points on $\surface$ is empty, despite there being adelic points on $\surface$ that are orthogonal to all elements in the Brauer group of $\surface$. A fortiori the Brauer group does not obstruct the existence of $\Q$-rational $0$-cycles of degree $1$ on $\surface$. In this short note we show that, as predicted by a conjecture of Colliot-Th\'el\`ene \cite{Colliot-Thelene99}, the surface $\surface$ does in fact possess a $\Q$-rational $0$-cycle of degree $1$.
	
	To wit, 
	\begin{align*}
		x_0 &= \theta^2 +1\\
		y_0 &= 3357\theta^2 - 2133\theta + 4851\\
		z_0 &= 2826\theta^2 - 2025\theta + 4158\\
		t_0 &= -42\theta^2 + 24\theta - 54\;,
	\end{align*}
	where $\theta$ satisfies $\theta^3 + \theta + 1 = 0$ are the coordinates of a closed point of $\surface$ whose residue field is the cubic number field $L = \Q[\theta]$. As there are obviously $0$-cycles of degree $4$ on $\surface$, this shows that there is a $0$-cycle of degree $1$ on $\surface$. 

	We note, however, that the far more important question of whether the Brauer-Manin obstruction to the existence of $0$-cycles of degree $1$ is the only one for bielliptic surfaces remains open. 

	Let us briefly explain how this point was discovered. Consider the genus one curves
	\[
		\curveC : U^2 = g(T) = 3(T^4-54T^2 - 117T-243)\,,
	\]
	and
	\[
		\curveD : \begin{cases} Y^2 = p(X) = X^2 + 1  &\;\\
					Z^2 = q(X) = X^2 + 2\,.
			\end{cases}
	\]
	There are actions of $\mu_2 = \{\pm1\}$ on $\curveC$ and $\curveD$ given by
	\[
		(T,U) \mapsto (T,-U) \quad \text{and} \quad (X,Y,Z) \mapsto (X,-Y,-Z)\,.
	\] 
	The quotient of $\curveC\times \curveD$ by the diagonal action of $\mu_2$ is the $\mu_2$-torsor $\rho:\curveC\times \curveD \to \surface$ given by $s^2 = g(t)$ over $\surface$. After a change of variables, $\surface$ may be written
	\[
		\surface : \begin{cases} y^2 = g(t)p(x) & \;\\ z^2 = g(t)q(x) &\;\end{cases}
	\]
	and $\rho$ is defined by $x = X, t = T, y = UY, z = UZ$.

	\subsection*{Lemma}
		There is a $0$-cycle of degree $1$ on $\surface$ if and only if there is an odd degree number field $K$ and $a \in K^\times$ such that genus one curves over $K$ given by 
		 \[
			\curveC^a : aU^2 = g(T)  \quad \text{and} \quad\quad 
			\curveD^a : \begin{cases} aY^2 = p(X) &\;\\
					aZ^2 = q(X) &\;
			\end{cases}
		\]
		both possess a $K$-rational point.
	
	\begin{proof}
		It is clear that there is a $0$-cycle of degree $1$ on $\surface$ if and only if $\surface(K) \ne \emptyset$ for some odd degree number field $K$. Let $\rho_K : (\curveC\times \curveD)_K \to \surface_K$ denote the base change to $K$. By descent theory every element of $\surface(K) = \surface_K(K)$ lifts to a $K$-rational point on some twist of the $\mu_2$-torsor $\rho_K : (\curveC\times \curveD)_K \to \surface_K$ by a cocycle in $H^1(K,\mu_2) = K^\times/K^{\times 2}$. A straightforward computation shows that the twist of $(\curveC\times \curveD)_K$ by $a \in K^\times/K^{\times 2}$ is the product of the curves in the statement of the proposition.
	\end{proof}

	The quotient of $\curveD$ by $\mu_2$ is the $\mu_2$-torsor $\phi: \curveD \to \curveD'$ given by $U^2 = p(X)$ over the genus one curve $\curveD'$ given by $W^2 = p(X)q(X)$. This curve has two rational points at infinity. Fixing either as the identity endows $\curveD'$ with the structure of an elliptic curve such that the other point is $2$-torsion. There are rational points on $\curveD$ above both of these, and so we may view $\phi:\curveD \to \curveD'$ as an isogeny of elliptic curves. The twists $\curveD^a$ which have a $K$-rational point correspond to the image of $\curveD'(K)$ under the connecting homomorphism $\delta_\phi$ in the exact sequence
	\[
		1 \to \mu_2 \to \curveD(K) \stackrel{\phi}\to \curveD'(K) \stackrel{\delta_\phi} \to H^1(K,\mu_2) \simeq K^\times/K^{\times 2}\,.
	\]
	In particular, for a given $K$, one can compute by means of an explicit $2$-isogeny descent a finite set $A \subset K^\times/K^{\times 2}$ such that every twist with a $K$-point is $\curveD^a$ for some $a \in A$. Thus the determination of $\surface(K)$ is reduced to the determination of the set of $K$-rational points on an explicit finite set of genus $1$ curves over $K$, namely the curves $\curveC^a$ and $\curveD^a$ with $a \in A$. Conjecturally, this is a finite computation. For $K$ of small degree and discriminant it is often possible in practice.

	We carried out these computations for various number fields using the Magma Computational Algebra System \cite{Magma}. There are no points on $\surface$ defined over the cubic field $L_1$ of smallest absolute discriminant. The map $\phi:\curveD(L_1) \to \curveD'(L_1)$ is surjective and so the only twist with $L_1$-points is $\curveD  = \curveD^1$. But as $\curveC = \curveC^1$ has no points over any odd degree number field we see that $\surface(L_1) = \emptyset$. The point $(x_0,y_0,z_0,t_0) \in \surface(L)$ is defined over the cubic field of second smallest absolute discriminant. The group $\curveD'(L)$ has rank $1$ and is generated by the points with $X$-coordinate equal to $x_0 = \theta^2 + 1$. This corresponds to the twist by $a = 6\theta^2 - 4\theta +9 \in L^\times$ and one finds that $\curveC^a$ has a point with $T$-coordinate equal to $t_0 = -42\theta^2 + 24\theta -54$. In fact, $\curveC^a(L)$ and $\curveD^a(L)$ are both infinite, showing that $\surface(L)$ is Zariski dense in $\surface$.
	
	\subsection*{Acknowledgements}
		The author would like to thank J.-L. Colliot-Th\'el\`ene for helpful comments and suggestions.

%%%%%%%%%%%%%%%%%%%%%%%%%%%%
%%%%%%%%%%%%%%%%%%%%%%%%%%%%
%% BIBLIOGRAPHY %%%%%%%%%%%%%%%%%
%%%%%%%%%%%%%%%%%%%%%%%%%%%%
%%%%%%%%%%%%%%%%%%%%%%%%%%%%

\end{document}